\documentclass[12pt]{article}

\usepackage{amsmath,amssymb}
\begin{document}

\pagestyle{myheadings} \markright{MELLIN TRANSFORMS OF
$P$-ADIC WHITTAKER FUNCTIONS}

\title{Mellin transforms of $p$-adic Whittaker functions}
\author{Anton Deitmar}

\date{}
\maketitle

$$ $$

\tableofcontents

\def \1{{\bf 1}}
\def \a{{{\mathfrak a}}}
\def \ad{{\rm ad}}
\def \al{\alpha}
\def \ar{{\alpha_r}}
\def \A{{\mathbb A}}
\def \Ad{{\rm Ad}}
\def \Aut{{\rm Aut}}
\def \b{{{\mathfrak b}}}
\def \bs{\backslash}
\def \B{{\cal B}}
\def \c{{\mathfrak c}}
\def \cent{{\rm cent}}
\def \C{{\mathbb C}}
\def \CA{{\cal A}}
\def \CB{{\cal B}}
\def \CC{{\cal C}}
\def \CD{{\cal D}}
\def \CE{{\cal E}}
\def \CF{{\cal F}}
\def \CG{{\cal G}}
\def \CH{{\cal H}}
\def \CHC{{\cal HC}}
\def \CL{{\cal L}}
\def \CM{{\cal M}}
\def \CN{{\cal N}}
\def \CP{{\cal P}}
\def \CQ{{\cal Q}}
\def \CO{{\cal O}}
\def \CS{{\cal S}}
\def \CT{{\cal T}}
\def \CV{{\cal V}}
\def \CW{{\cal W}}
\def \det{{\rm det}}
\def \diag{{\rm diag}}
\def \dist{{\rm dist}}
\def \End{{\rm End}}
\def \eps{\varepsilon}
\def \eqn{\begin{eqnarray*}}
\def \endeqn{\end{eqnarray*}}
\def \F{{\mathbb F}}
\def \Fx{{\mathfrak x}}
\def \FX{{\mathfrak X}}
\def \g{{{\mathfrak g}}}
\def \ga{\gamma}
\def \Ga{\Gamma}
\def \GL{{\rm GL}}
\def \h{{{\mathfrak h}}}
\def \Hom{{\rm Hom}}
\def \im{{\rm im}}
\def \Im{{\rm Im}}
\def \Ind{{\rm Ind}}
\def \k{{{\mathfrak k}}}
\def \K{{\cal K}}
\def \l{{\mathfrak l}}
\def \la{\lambda}
\def \lap{\triangle}
\def \li{{\rm li}}
\def \Lie{{\rm Lie}}
\def \m{{{\mathfrak m}}}
\def \mod{{\rm mod}}
\def \n{{{\mathfrak n}}}
\def \name{\bf}
\def \Mat{{\rm Mat}}
\def \N{\mathbb N}
\def \o{{\mathfrak o}}
\def \ord{{\rm ord}}
\def \O{{\cal O}}
\def \p{{{\mathfrak p}}}
\def \ph{\varphi}
\def \prf{\noindent{\bf Proof: }}
\def \Per{{\rm Per}}
\def \PGL{{\rm PGL}}
\def \q{{\mathfrak q}}
\def \qed{\ifmmode\eqno Q.E.D.\else\noproof\vskip 12pt plus 3pt minus 9pt \fi}
 \def\noproof{{\unskip\nobreak\hfill\penalty50\hskip2em\hbox{}%
     \nobreak\hfill Q.E.D.\parfillskip=0pt%
     \finalhyphendemerits=0\par}}
\def \Q{\mathbb Q}
\def \res{{\rm res}}
\def \R{{\mathbb R}}
\def \Re{{\rm Re \hspace{1pt}}}
\def \r{{\mathfrak r}}
\def \ra{\rightarrow}
\def \rank{{\rm rank}}
\def \sign{{\rm sign}}
\def \supp{{\rm supp}}
\def \SL{{\rm SL}}
\def \Spin{{\rm Spin}}
\def \SU{{\rm SU}}
\def \t{{{\mathfrak t}}}
\def \T{{\mathbb T}}
\def \tr{{\hspace{1pt}\rm tr\hspace{2pt}}}
\def \vol{{\rm vol}}
\def \z{\zeta}
\def \Z{\mathbb Z}
\def \={\ =\ }

\newcommand{\frack}[2]{\genfrac{}{}{0pt}{}{#1}{#2}}
\newcommand{\rez}[1]{\frac{1}{#1}}
\newcommand{\der}[1]{\frac{\partial}{\partial #1}}
\newcommand{\norm}[1]{\parallel #1 \parallel}
\renewcommand{\matrix}[4]{\left( \begin{array}{cc}#1 & #2 \\ #3 & #4 \end{array}
            \right)}
\renewcommand{\sp}[2]{\langle #1,#2\rangle}

\newcounter{lemma}
\newcounter{corollary}
\newcounter{proposition}
\newcounter{theorex}

\newtheorem{conjecture}{\hspace{-20pt}\stepcounter{lemma} \stepcounter{corollary}
    \stepcounter{proposition}\stepcounter{theorex}Conjecture}[section]
\newtheorem{lemma}{\hspace{-20pt}\stepcounter{conjecture}\stepcounter{corollary}
    \stepcounter{proposition}\stepcounter{theorex}Lemma}[section]
\newtheorem{corollary}{\hspace{-20pt}\stepcounter{conjecture}\stepcounter{lemma}
    \stepcounter{proposition}\stepcounter{theorex}Corollary}[section]
\newtheorem{proposition}{\hspace{-20pt}\stepcounter{conjecture}\stepcounter{lemma}
    \stepcounter{corollary}\stepcounter{theorex}Proposition}[section]

\newtheorem{theorex}{\hspace{-40pt}\stepcounter{conjecture} \stepcounter{lemma}
    \stepcounter{corollary} \stepcounter{proposition}Theorem}[section]
\newenvironment{theorem}{\vspace{30pt}\begin{theorex}}{\end{theorex}\vspace{30pt}}

\newpage
\begin{center}
{\bf Introduction}
\end{center}
Whittaker functions occur naturally in the theory of
automorphic forms as the Fourier coefficients of cusp
forms. As a consequence, inner products of Poincar\'e
series with cusp forms are Mellin transforms of Whittaker
functions \cite{bump-friedberg-goldfeld,selberg}. Mellin
transforms of Whittaker functions have been considered
mainly in the context of automorphic $L$-functions for the
group $\GL_n$, see \cite{Bump-Friedberg,Friedberg-Goldfeld,
jacquet, Jacquet-PS-Shalika, JL,Stade}. Some authors
consider Mellin transforms with respect to the central
torus and some restrict to the archimedean case. In
\cite{Friedberg-Goldfeld} the archimedean transform in the
unramified case is considered. In \cite{Deitmar:Whitt} the
author established meromorphicity of the Mellin transform
for arbitrary archimedean groups and arbitrary
representations.

In this paper the nonarchimedean situation is considered.
At first it has to be secured that the Mellin transforms
actually do exist by showing that Whittaker vectors are of
moderate growth (Theorem \ref{2.2}). Next it is shown
that, for a dense set of vectors, the Mellin transform is
a rational function in $q$-powers, where $q$ is the order
of the residue class field. Finally, for unramified
representations, the class one Mellin transform is
computed explicitly, showing that it is a product of a
simple entire function and an Euler factor attached to the
adjoint representation.

The motivation for this paper is to seek an analytic
approach to general automorphic $L$-functions. Such an
approach would require a representation of a given
$L$-function as, say, a Mellin integral. Up to now this
essentially only can be done for standard $L$-functions of
the group $\GL_n$. The Mellin transform introduced here
has the advantage that it is defined for any reductive
group.

\newpage
\section{Whittaker functions}
Let $F$ be a nonarchimedean local field with valuation
ring $\CO$ and uniformizer $\varpi$. Denote by $G$ be a
connected reductive group over $F$. Let $K\subset G$ be a
good maximal compact subgroup. Choose a minimal parabolic
subgroup $P=MN$ of $G$ with Levi component $M$. Let $A$
denote the greatest split torus in the center of $M$. Then
$A$ is called the {\it split component} of $P$. Let
$\Phi=\Phi(G,A)$ be the root system of the pair $(G,A)$,
i.e. $\Phi$ consists of all homomorphisms $\alpha :A\ra
\GL_1$ such that there is $X$ in the Lie algebra of $G$
with $\Ad(a)X= a^\alpha X$ for every $a\in A$. Given
$\alpha$, let $\n_\alpha$ be the Lie algebra generated by
all such $X$ and let $N_\alpha$ be the closed subgroup of
$N$ corresponding to $\n_\alpha$. Let $\Phi^+=\Phi(P,A)$ be
the subset of $\Phi$ consisting of all positive roots with
respect to $P$. Let $\Delta\subset\Phi^+$ be the subset of
simple roots. A character $\eta : N\ra\T$ is called {\it
generic} if $\eta$ is nontrivial on $N_\alpha$ for every
$\alpha\in\Delta$. Let $A^-\subset A$ be the set of all
$a\in A$ such that $|a^\alpha| <1$ for any $\alpha\in
\Delta$.

Denote by $N(A)$ be the normalizer of $A$ in $G$. The {\it
Weyl group}:
$$
W\= N(A)/A
$$
is a finite group which acts on $A$. Each element $w$ of $W$ has
a representative in $K$. We frequently use the same notation for
the element $w$ of $W$ and its representative in $K$.

Let $X^*(A)=\Hom(A,\GL_1)$ be the group of all
homomorphisms as algebraic groups from $A$ to $\GL_1$. This
group is isomorphic to $\Z^r$ with $r=\dim A$. Likewise let
$X_*(A)=\Hom(\GL_1,A)$. There is a natural $\Z$-valued
pairing
\eqn
X^*(A)\times X_*(A) &\ra& \Hom(\GL_1,\GL_1)\cong\Z\\
(\alpha,\eta) &\mapsto& \alpha\circ\eta.
\endeqn
For every root $\alpha\in\Phi(A,G)\subset X^*(A)$ let
$\breve{\alpha}\in X_*(A)$ be its coroot. Then
$(\alpha,\breve{\alpha})=2$. The valuation $v$ of $F$
gives a group homomorphism $\GL_1(F)\ra\Z$. Let $A_c$ be
the unique maximal compact subgroup of $A$. Let $\Sigma
=A/A_c$; then $\Sigma$ is a $\Z$-lattice of rank $r=\dim
A$. By composing with the valuation $v$ the group $X^*(A)$
can be identified with
$$
\Sigma^*\=\Hom(\Sigma,\Z).
$$
Let
$$
\a_0^*\=\Hom(\Sigma,\R)\ \cong\ X^*(A)\otimes\R
$$
be the real vector space of all group homomorphisms from
$\Sigma$ to $\R$ and let
$\a^*=\a_0^*\otimes\C=\Hom(\Sigma,\C)\cong
X^*(A)\otimes\C$. For $a\in A$ and $\la\in\a^*$ let
$$
a^\la\=q^{-\la(a)},
$$
where $q$ is the number of elements in the residue class field of
$F$. In this way we get an identification
$$
{\a^*}/\mbox{\small$\frac{2\pi i}{\log q}$}\Sigma^* \
\cong\ \Hom(\Sigma,\C^\times).
$$
A quasicharacter $\nu : A\ra\C^\times$ is called {\it
unramified} if $\nu$ is trivial on $A_c$. The set
$\Hom(\Sigma,\C^\times)$ can be identified with the set of
unramified quasicharacters on $A$. Any unramified
quasicharacter $\nu$ can thus be given a unique real part
$\Re(\nu)\in \a_0^*$. This definition extends to not
necessarily unramified quasicharacters $\chi:A\ra\C^\times$
as follows. Choose a splitting $s:\Sigma\ra A$ of the exact
sequence
$$
1\ra A_c\ra A\ra\Sigma\ra 1.
$$
Then $\nu=\chi\circ s$ is an unramified character of $A$. Set
$$
\Re(\chi)\=\Re(\nu).
$$
This definition does not depend on the choice of the splitting
$s$. For quasicharacters $\chi$, $\chi'$ and $a\in A$ we will
frequently write $a^\chi$ instead of $\chi(a)$ and
$a^{\chi+\chi'}$ instead of $\chi(a)\chi'(a)$. Note that the
absolute value satisfies $|a^\chi|=a^{\Re(\chi)}$ and that a
quasicharacter $\chi$ actually is a character if and only if
$\Re(\chi)=0$.

Let $\Delta_P : P\ra\R_+$ be the modular function of the
group $P$. Then there is $\rho\in\a_0^*$ such that
$\Delta_P(a)=a^{2\rho}$. For $\nu\in\a^*$ and a root
$\alpha$ let
$$
\nu_\alpha\= (\nu,\breve{\alpha})\ \in\
X^*(\GL_1)\otimes\C\ \cong\ \C.
$$
Note that $\nu\in\a_0^*$ implies $\nu_\alpha\in\R$ for
every $\alpha$. For $\nu\in\a_0^*$ we say that $\nu$ is
positive, $\nu>0$, if $\nu_\alpha>0$ for every positive
root $\alpha$.

{\bf Example.} Let $G=\GL_n(F)$ and let $\varpi_j\in G$ be
the diagonal matrix
$\varpi_j=\diag(1,\dots,1,\varpi,1,\dots,1)$ with the
$\varpi$ on the $j$-th position. Let $\nu\in\a^*$ and let
$$
\nu_j\=\nu(\varpi_j A_c)\ \in\ \C.
$$
Let $\alpha$ be a root, say
$\alpha(\diag(a_1,\dots,a_n))=\frac{a_i}{a_j}$. Then
$$
\nu_\alpha\= \nu_i-\nu_j.
$$
Hence $\nu\in\a_0^*$ is positive if and only if
$\nu_1>\nu_2>\dots >\nu_n$.

By a {\it representation} of $G$ we will simply mean a group
homomorphism $\pi:G\ra GL(V)$ from $G$ to the group of all linear
homomorphisms of a complex vector space $V$, in other words, $G$
acts linearly on $V$. The representation $\pi$ is called {\it
irreducible} if there is no proper subrepresentation; i.e. if the
only $G$-stable subspaces of $V$ are the zero space and $V$
itself. The representation $\rho$ is called {\it smooth} if every
vector $v\in V$ is fixed by some open subgroup $H$ of $G$. The
representation $\pi$ is called {\it admissible} if it is smooth
and for every open subgroup $H$ of $G$ the space of fixed vectors
$V^H$ is finite dimensional. Finally a representation $\pi$ is
called {\it quasi-simple} if the center of $G$ acts by scalar
multiples of the identity.

Let $(\pi,V_\pi)$ be an irreducible admissible representation of
$G$. Let $\eta :N\ra\T=\{ z\in \C\ :\ |z|=1\}$ be a character. A
nonzero linear functional $\psi=\psi_\eta : V_\pi\ra\C$ on $\pi$
is called a {\it Whittaker functional} with respect to $\eta$ if
for any  $n\in N$ and $v\in V_{\pi}$ we have
$\psi(\pi(n)v)\=\eta(n)\psi(v)$. Define the corresponding {\it
Whittaker function} by
$$
W_v(x)\=\psi(\pi(x)v),\ \ x\in G.
$$
The character $\eta$ is called {\it generic} if for every
$\alpha\in\Delta$ the restriction of $\eta$ to the group
$N_\alpha$ is nontrivial.

\begin{lemma}\label{1.1}
Suppose that $\eta$ is generic. Then, for any $v\in V_\pi$ there
is an element $a_0\in A$ such that the Whittaker function $W_v|_A$
has support in $a_0A^-$.
\end{lemma}

\prf Let $v\in V_\pi$. Since the stabilizer of $v$ in $G$ is
open, there is a compact open subgroup $N_0$ of $N$ such that for
any $n\in N_0$ we have $\pi(n)v=v$. This implies that for $a\in
A$ and $n\in N_0$:
\eqn
W_v(a) &=& \psi(\pi(a)v)\\
&=& \psi(\pi(a)\pi(n)v)\\
&=& \psi(\pi(^{a}n)\pi(a)v)\\
&=& \eta(^{a}n)\psi(\pi(a)v)\\
&=& \eta(^{a}n)W_v(a),
\endeqn
where we have written $^{a}n$ for $na^{-1}$. So $W_v(a)\ne
0$ implies $\eta(^{a}n)=1$ for every $n\in N_0$. The
assertion follows.
\qed

In \cite{CassShal} it is shown that if $\pi$ is a supercuspidal
representation then $W_v(a)$ has compact support modulo the
center of $G$.

\section{Induced representations}
Let $\theta$ be a subset of $\Delta$, and let
$P_\theta=M_\theta N_\theta$ be the corresponding
parabolic containing $P$. Let $A_\theta$ be the split
component of $P_\theta$. The extremal cases are
$P_\emptyset =P$ and $P_\Delta =G$. Note that
$A_\Delta=A_G$ is the largest split torus in the center of
$G$. Let $A^{der}=A\cap G^{der}$, where $G^{der}$ is the
derived group of $G$. The product $A_G A^{der}$ has finite
index in $A$. Let $x\mapsto x^{2\rho_{_\theta}}$ be the
modular function of $P_\theta$. Let $\sigma$ be an
admissible representation of $M_\theta$. Denote by
$\pi=\pi_\sigma$ the representation smoothly induced from
$\sigma\otimes 1$ on $P_\theta$. The space of this
representation is the space $V_\pi$ of all locally
constant functions $f:G\ra V_\sigma$ which satisfy $f(mnx)
= m^{\rho_{_\theta}}\sigma(m)f(x)$ for all $mn\in P_\theta$
and all $x\in G$. The representation $\pi$ then is defined
by $\pi(x)f(y)=f(yx)$. Recall that if $\sigma$ is
quasi-simple of finite length, then $\pi_\sigma$ is a
quasi-simple representation of finite length.

By a theorem of Harish-Chandra every irreducible admissible
representation $\pi$ is a quotient of some induced representation
$\pi_\sigma$, for $\sigma$ supercuspidal. So any Whittaker
function for $\pi$ will induce one for $\pi_\sigma$. For the
analysis of Whittaker functions it is therefore sufficient only
to consider induced representations.

Let $\bar P_\theta=M_\theta \bar N_\theta$ be the parabolic
opposite to $P_\theta$. This means that $\bar N_\theta$ is given
as follows. Let $\Phi^{+}_\theta$ be the intersection of $\Phi^+$
with the linear span of $\theta$. Then $N_\theta
=\prod_{\alpha\in\Delta\setminus\Phi_\theta^{+}} N_\alpha$. Now
define
$$
\bar N_\theta\=\prod_{\alpha\in\Delta\setminus\Phi_\theta^{+}}
N_{-\alpha}.
$$
There is a different way to construct the group $\bar N_\theta$
which is quite useful. Let $w_{\theta,l}$ denote the longest
element of the Weyl group $W_\theta\subset W$ of $A_\theta$. Then
for $\theta=\emptyset$ we get $w_l=w_{\emptyset,l}$, which is the
longest element of $W$. Set $w_\theta=w_{\theta,l}w_l$. and set
$\bar\theta=w_\theta^{-1}\theta$, then it turns out
\cite{CassShal} that $\bar N_\theta =w_\theta
N_{\bar\theta}w_\theta^{-1}$, and that $w_\theta$ is the longest
element of the set $[W_\theta\bs W]=\{ w\in W | w^{-1}\theta
>0\}$.

Let $\pi=\pi_\sigma$, where $\sigma$ is not necessarily
supercuspidal but quasisimple and of finite length. Let
$V_\pi^c\subset V_\pi$ be the subspace of functions $f:G\ra
V_\sigma$ such that $n\mapsto f(w_\theta n)$ has compact support
in $N_{\bar\theta}$. Note that this subspace is stable under the
action of $P_\theta$, but not in general under the action of the
group $G$. Let $\psi_\sigma$ be a Whittaker functional on
$V_\sigma$ for the character $^{w_\theta}\eta$, i.e.
$\psi_\sigma$ satisfies
$\psi_\sigma(\sigma(n)v)=\eta(w_\theta^{-1}
nw_\theta)\psi_\sigma(v)$ for every $n\in N\cap M_\theta$. Then
the linear map $\psi:V_\pi^c\ra\C$ given by
$$
\psi(f)\=\int_{N_{\bar\theta}} \eta(n^{-1})\psi_\sigma(f(w_\theta
n) dn
$$
satisfies $\psi(\pi(n)f)=\eta(n)\psi(f)$ for any $n\in N$.

\begin{theorem}\label{rodier}
(Rodier) The functional $\psi$ extends to a unique Whittaker
functional on $V_\pi$. The map sending $\psi_\sigma$ to $\psi$ is
a bijection from the space of all Whittaker functionals on
$V_\sigma$ to the space of all Whittaker functionals on $V_\pi$.
\end{theorem}

\prf The assertion follows from Theorem 1.6 of \cite{CassShal}.
\qed

\section{Moderate growth}

We fix an embedding $G\hookrightarrow \GL_n(F)$ with
$K\subset \GL_n(\CO)$. On $Mat_n(F)\supset \GL_n(F)$ we
have a norm
$$
\norm{A}\=\max_{i,j} |a_{i,j}|,
$$
which satisfies $\norm{AB}\le\norm A\norm B$ for all $A,B\in
Mat_n(F)$.

We say that a function $f$ on $G$ is {\it of moderate
growth} if there are $C>0$ and $m\in\N$ such that for
every $x\in G$ we have
$$
|f(x)|\ \le\ C\norm x^m.
$$

\begin{theorem}\label{2.2}
Let $\eta$ be generic, and let $\pi$ be an irreducible admissible
representation of $G$. Then for any $f\in V_\pi$ the Whittaker
function $W_f$ is of moderate growth.
\end{theorem}

\prf According to a theorem of Harish-Chandra, every irreducible
admissible representation is a quotient of a representation
parabolically induced from a supercuspidal representation. It
therefore suffices to prove the claim not for irreducible
representations but for induced ones.

We have the decomposition $G=PK=NMK$. Let $f\in V_\pi$ as
above; then
$W_f(nmk)=\psi(\pi(nmk)f)=\eta(n)\psi(\pi(m)\pi(k)f)$.
Since $K$ is compact, it suffices to show that $W_f$ is of
moderate growth on $M$. Since further $M=A\omega$ for some
compact set $\omega$ it follows that we only have to check
that $W_f$ is of moderate growth on $A$.

The proof of the theorem will be based on induction. At first
assume that $\pi$ is supercuspidal, then $W_f(a)$ is of compact
support, which implies the claim in this case.

Next assume that $\pi$ is not supercuspidal. Then there is
a supercuspidal representation $\sigma$ of some
Levi-component $M_{\theta'}$ such that $\pi$ is a quotient
of $\pi_\sigma=\Ind_{P_{\theta'}}^G(\sigma\otimes 1)$.
Thus it suffices to show the claim in the case
$\pi=\pi_\sigma$. In this case $\pi$ is not in general
irreducible, but still quasi-simple and of finite length.
Let $P_\theta$ be a maximal proper parabolic with
$P_{\theta'}\subset P_\theta$. Let $\tau$ be the
representation of $M_\theta$ induced from $\sigma\otimes
1$ on $P_{\theta'}\cap M_\theta$. By the double induction
formula we have that $\pi=\pi_\tau$, the representation
induced from $\tau$. By induction hypothesis we may assume
that any Whittaker functional for $\tau$ is of moderate
growth. So let $\pi=\pi_\tau$, then $V_\pi$ is the space
of all locally constant functions $f:G\ra V_\tau$ such
that $f(mnx)=m^{\rho_\theta}\tau(m) f(x)$ for $mn\in
P_\theta$.

\begin{lemma}
Let $N_c\subset N$ be a compact open subgroup such that $\eta$ is
nontrivial on $N_c\cap N_\alpha$ for every $\alpha\in\Delta$. For
$f\in V_\pi$ define
$$
\CP(f)(x)\=\rez{\vol(N_c)}\int_{N_c}\eta(n^{-1})f(xn)dn.
$$
Then $\CP(f)$ has support in the open Bruhat cell $P_\theta
w_\theta P$.
\end{lemma}

\prf According to the Bruhat decomposition we have $G=P_\theta
w_\theta P\cup P_\theta$, since $P_\theta$ is maximal. So we have
to show that $\CP(f)(1)=0$. For this recall
$$
\CP(f)(1)\=\rez{\vol(N_c)}\int_{N_c}\eta(n^{-1})f(n) dn.
$$
Now let $n_1\in N_c\cap N_\theta$ such that $\eta(n_1)\ne 1$.
Since $f(n)=f(n_1n)$ it follows that
$\CO(f)(1)=\eta(n_1)\CP(f)(1)$, so $\CP(f)(1)=0$. The lemma
follows.
\qed

For $j=0,1,2,\dots$ let $\Ga(j)$ denote the subgroup of
$\GL_n(\CO)$ consisting of all elements congruent to the
unit matrix modulo $\varpi^{2j}\CO$. Let $K_j=K\cap\Ga(j)$
then $K_j$ is normal in $K$ and the sequence $(K_j)_j$
forms a neighbourhood basis of the unit in $G$. Let
$\log_q(x)=\frac{\log(x)}{\log(q)}$. Then for $a\in A$ the
number
$$
e(a)\= \max_{w\in W}\log_q(a^{2w\rho})
$$
is in $\N_0$. A function $\ph$ on $A$ is of moderate growth if
and only if there are $\mu\in\N$ and $C>0$ such that $|\ph(a)|\le
Cq^{\mu e(a)}$. Note that there is a natural number $k_0$ such
that for any $a\in A$ we have $aK_ja^{-1}\supset K_{j+k_0 e(a)}$.
Let $V_\pi(j)= V_\pi(K_j)$ be the set of vectors in $V_\pi$ which
are fixed by $K_j$. Then it follows that for each $a\in A$ we
have $\pi(a)V_\pi(j)\subset V_\pi(j+k_0 e(a))$.

\begin{lemma}\label{2.4}
There is a natural number $\tilde m$ such that for each $j\in\N$
the set $P_\theta w_\theta N\setminus P_\theta K_j$ is a subset of
$P_\theta w_\theta N_{\tilde m j}$.
\end{lemma}

\prf Let $\bar P=M\bar N$ be the opposite parabolic and let $\bar
N_0=\bar N\cap K$. Let $a\in A^+$ be so large that $a\bar N_0
a^{-1}\subset K_j$. Then it follows that
\eqn
P_\theta w_\theta N\setminus P_\theta K_j &\subset&
P_\theta w_\theta N\setminus P_\theta  a\bar N_0 a^{-1}\\
&=& [P_\theta wN\setminus P_\theta \bar N_0]
a^{-1}.
\endeqn
So, if $ P_\theta wN\setminus P_\theta \bar N_0 $ is
contained in, say, $P_\theta wN_j$, then $ P_\theta
wN\setminus P_\theta K_j
$
is contained in $P_\theta w(aN_j a^{-1})$. The lemma follows.
\qed

\begin{lemma}\label{2.7}
There is a natural number $k_1$ such that $\CP(V_\pi(j))\subset
V_\pi(j+k_1)$.
\end{lemma}

\prf We have
$$
\CP(V_\pi(j))\= \CP(V_\pi(K_j)\ \subset\ V_\pi\left(
\bigcap_{n\in N_c}nK_jn^{-1}\right).
$$
Let $a_0\in A$ be such that $N_c\subset a_0 N_0a^{-1}$. then
\eqn
\bigcap_{n\in N_c}nK_jn^{-1} &\supset &\bigcap_{n_0\in N_0}
a_0n_0a_0^{-1} K_j a_0n_0^{-1}a_0^{-1}\\
&\supset& \bigcap_{n_0\in N_0}
a_0n_0 K_{j+k_0e(a_0)} n_0^{-1}a_0^{-1}.
\endeqn
Now $N_0$ is a subset of $K$, and $K_j$ is normal in $K$.
Therefore the latter equals
$
a_0K_{j+k_0 e(a_0)}a_0^{-1}\ \supset\ K_{j+2k_0e(a_0)}.
$
With $k_1=2k_0e(a_0)$ the lemma follows.
\qed

Now we continue the proof of the theorem. For every $f\in V_\pi$
we have
\eqn
\psi(\CP(f)) &=&
\rez{\vol(N_c)}\int_{N_c}\eta(n^{-1})\psi(\pi(n)f)dn\\
&=& \psi(f)\rez{\vol(N_c)}\int_{N_c}dn\\
&=& \psi(f).
\endeqn
Suppose $f\in V_\pi(j)$; then for $a\in A$ we have that
$\pi(a)f\in V_\pi(j+k_0 e(a))$, and by Lemma \ref{2.7} it follows
that
$
\CP(\pi(a)f)\ \in\ V_\pi(j+k_0 e(a)+k_1).
$
By Lemma \ref{2.4} this implies that
$
\supp(\CP(\pi(a)f))\ \subset\ P_\theta w_\theta N_{\tilde
m(j+k_0e(a)+k_1)}.
$
Therefore
\eqn
W_f(a) &=& \psi(\pi(a)f)\\
&=& \psi(\CP(\pi(a)f))\\
&=& \int_{N_{\bar\theta}\cap N_{\tilde m(j+k_0e(a)+k_1)}}
\eta(n^{-1})\psi_\tau(f(w_\theta n a)) d n,
\endeqn
for some Whittaker functional $\psi_\tau$ on $V_\tau$. This last
line is implied by Theorem \ref{rodier}.

For an element $x$ of $G$ we write the $MNK$-decomposition as
$$
x=\underline m(x)\underline n(x)\underline k(x).
$$
Then
$$
\psi_\tau(f(w_\theta na))\= \underline m(w_\theta
na)^{\rho_\theta}\psi_\tau(\tau(\underline m(w_\theta
na))f(\underline k(w_\theta na))).
$$
The function $na\mapsto f(\underline k(w_\theta na))$ takes only
finitely many values since $f$ is locally constant. By induction
hypothesis the Whittaker function $m\mapsto \psi_\tau(\tau(m) v)$
is of moderate growth for any $v\in V_\tau$. It follows that
there is $m_1\in\N$ and $C>0$ such that
$$
\psi_\tau(f(w_\theta na))\ \le\ C\norm{\underline m(w_\theta
na)}^{m_1}.
$$
For any $x\in \GL_n(F)$ and $k_1,k_2\in \GL_n(\CO)$ we have
$$
\norm{k_1 xk_2}\ \le\ \norm{k_1}\norm x\norm{k_2}\=\norm x.
$$
Varying $x,k_1,k_2$ we get
$$
\norm{k_1 xk_2}\=\norm x.
$$
There is a parabolic $P'_\theta$ of $\GL_n(F)$ such that
$P_\theta=P'_\theta\cap G$. Modulo conjugation by
$\GL_n(\CO)$ we may assume that $P'_\theta$ is a standard
parabolic, i.e. it contains the upper triangular matrices.
We can choose the Levi component $M'_\theta$ of
$P'_\theta$ with $M_\theta=M'_\theta\cap G$ as a standard
component. From this it follows that for any $m\in
M_\theta$ and $n\in N_\theta$ we have
$
\norm m\ \le\ \norm{mn}.
$
In particular
$
\norm{\underline m(w_\theta na)}\ \le\ \norm{\underline
m(w_\theta na)\underline n(w_\theta na)},
$
which by the above equals
$
\norm{w_\theta na}\=\norm{na}.
$
It follows that
\eqn
\psi_\tau(f(w_\theta na)) &\le& C\norm{na}^{m_1}\\
&\le& C\norm n^{m_1} \norm a^{m_1}.
\endeqn
If $n\in N_{\tilde m(j+k_0 e(a) +k_1)}$, then
$
\norm n \ \le\ q^{\tilde m(j+k_0 e(a) +k_1)}.
$
Therefore
\eqn
|W_f(a)| &\le& C\norm a^{m_1} \int_{N_{\bar\theta}\cap N_{\tilde
m(j+k_0e(a)+k_1)}} \norm n^{m_1} d n\\
&\le& C\norm a^{m_1} q^{m_1\tilde
m(j+k_0e(a)+k_1)}\vol(N_{\bar\theta}\cap N_{\tilde
m(j+k_0e(a)+k_1)})
\endeqn
There is $m_2\in\N$ such that
$$
\vol(N_{\bar\theta}\cap N_{\tilde m(j+(k_0)e(a)+k_1)})\ \le\
q^{m_2\tilde m(j+(k_0)e(a)+k_1)}.
$$
So that finally
\eqn
|W_f(a)|&\le& C\norm a^{m_1}q^{(m_1+m_2)\tilde
m(j+(k_0)e(a)+k_1)}\\
&\le& C_1\norm a^{m_1+\tilde m k_0 m_3},
\endeqn
for some $m_3\in\N$ and $C_1>0$.
\qed

\section{The Mellin transform}
For every simple root $\alpha$ let $A_\alpha$ denote the
set of $a\in A^{der}$ with $a^\beta=1$ for every
$\beta\in\Delta$, $\beta\ne \alpha$. Then the product map
induces an exact sequence
$$
1\ra E\ra \prod_{\alpha\in\Delta} A_\alpha\ra A^{der}\ra
F\ra 1,
$$
where $E$ and $F$ are finite abelian groups. Let $A'$ be
the image of $\prod_{\alpha\in\Delta}A_\alpha$ in
$A^{der}.$ The root lattice $R$ is in general a sublattice
of the character lattice $X^*(A^{der})$. The restriction
map however gives rise to an isomorphism $R\cong X^*(A')$.
By this reason the torus $A'$ will be called the {\it root
torus}.

Let $\eta$ be generic, and let $\pi$ be an irreducible
admissible representation of $G$. Let $\la_0 : A'\ra\T$ be
a character. For $\la\in\a^*$ and $f\in V_\pi$ set
$$
I_{\la_0,f}(\la)\=\int_{A'} W_f(a) a^{\la_0+\la-\rho} da.
$$
From Theorem \ref{2.2} and Lemma \ref{1.1} it follows that
the integral $I_f(\la)$ converges locally uniformly for
$\Re(\la)>>0$.\\
If $\la_0$ is the trivial character we also write
$I_f(\la)$.

\begin{theorem}\label{main}
Suppose that $\eta$ is generic and let $\pi=\pi_\sigma$,
where $\sigma$ is supercuspidal. Then for any $f\in
V_\pi^c$ the function $\la\mapsto I_{\la_0,f}(\la)$ extends
to a meromorphic function on $\a^*$. Indeed it is a
rational function in $(q^{\la_\alpha})_{\alpha\in\Delta}$.
\end{theorem}

\prf The function $f$ is compactly supported in $w_\theta
N_\theta$, and
\eqn
W_f(a)&=& \int_{N_{\bar\theta}}\eta(n^{-1})\psi_\sigma(f(w_\theta n a) dn\\
&=&
a^{2\rho_{_{\theta}}}\int_{N_{\bar\theta}}\eta(^an^{-1})\psi_\sigma(\sigma(^wa)f(w_\theta
n)) d n.
\endeqn
So that
$$
I_{\la_0,f}(\la)\= \int_{A'}
a^{\la_0+\la-\rho+2\rho_\theta}\int_{N_{\bar\theta}}\eta(^an^{-1})
\psi_\sigma(\sigma(^{w_\theta}a)f(w_\theta n))dnda.
$$
Note that the map $m\mapsto\sigma(^{w_\theta} m)$ is a
representation of $w_\theta^{-1}M_\theta w_\theta
=M_{\bar\theta}$.

Let $A_{CM}$ be the image of
$\prod_{\alpha\in\Delta\setminus\bar\theta}A_\alpha$ in
$A^{der}$. Then $A_{CM}$ lies in the center of
$M_{\bar\theta}$. Likewise let $A_M$ be the image of
$\prod_{\alpha\in\bar\theta}A_\alpha$. Then $A_M$ has a
finite intersection with the center of $M_{\bar\theta}$.
Write any $a\in A'$ as a product $a=a_{CM}a_M$, where
$a_{CM}\in A_{CM}$ and $a_M\in A_M$. Then
$$
\psi_\sigma(\sigma(^{w_\theta}a) f(w_\theta n))\=
a_{CM}^{w_\theta\chi}\psi_\sigma(\sigma(^{w_\theta}a_M)
f(w_\theta n))
$$
for a character $\chi=\chi_\sigma$, the central character of
$\sigma$. The locally constant function
$$
(a_M,n)\mapsto \psi_\sigma(\sigma(^{w_\theta}a_M) f(w_\theta n))
$$
has compact support in $A_M\times N_{\bar \theta}$. It therefore
is a linear combination of functions of the type $\1_{a_0
C}\times\1_{n_0V}$, where $a_0\in A_M$, $n_0\in N_{\bar\theta}$;
$C$ and $V$ are compact neighbourhoods of the units in $A_M$ and
$N_{\bar\theta}$. So we have to worry about expressions of the
form
$$
\int_{A_M}\int_{A_{CM}} (a_{CM} a_M)
^{w_\theta\chi+\la_0+\la-\rho+2\rho_\theta}
 \hspace{150pt} $$ $$ \hspace{50pt}\times
\int_{N_{\bar\theta}}\eta(^{a_Ma_{CM}}n^{-1})\1_{a_0
C}(a_M)\1_{n_0V}(n)dnda_M da_{CM}.
$$
For $C$ sufficiently small this equals
$$
a_0^{w_\theta\chi+\la_0+\la-\rho+\rho_\theta} \int_C
a_M^{w_\theta\chi+\la_0+\la-\rho+2\rho_\theta}
da_M\hspace{100pt}
$$ $$
\hspace{40pt}\times \int_{A_{CM}}
a_{CM}^{w_\theta\chi+\la_0+\la-\rho+2\rho_\theta}\eta(^{a_{CM}
a_0}n_0) \int_{V}\eta(^{a_{CM}}n^{-1})dn da_{CM}.
$$
The integral over the compact set $C$ is a polynomial in the
$q^{\pm\la_\alpha}$. If the representation $\pi$ itself is
supercuspidal then the last factor does not occur and the proof
is finished. So, let's assume that $\pi$ is not supercuspidal.
Then we have to treat the factor
$$
\int_{A_{CM}}
a_{CM}^{w_\theta\chi+\la_0+\la-\rho+2\rho_\theta}\eta(^{a_{CM}
a_0}n_0) \int_{V}\eta(^{a_{CM}}n^{-1})dn da_{CM}.
$$
Let $\Phi_{\bar\theta}$ be the set of all roots $\alpha\in\Phi$
which lie in the $\Z$-span of $\bar\theta$ and let
$\Phi_{\bar\theta}^+=\Phi_{\bar\theta}\cap\Phi^+$. Let
$N_{\bar\theta,0}$ be the subgroup of $N_{\bar\theta}$, generated
by all $N_\alpha$ where $\alpha\in\Phi_\theta^+$, but
$\alpha\notin \N\theta$. Then, as a set, we can write
$$
N_{\bar\theta}\= N_{\bar\theta,0}\times
\prod_{\alpha\in\Delta\setminus\bar\theta}N_\alpha.
$$
The Haar-measure of $N_{\bar\theta}$ can also be written as a
product of the Haar-measures of the subgroups. Now note that
$N_{\bar\theta,0}$ is a subgroup of the commutator group $[N,N]$
and so the character $\eta$ is trivial on $N_{\bar\theta,0}$. We
perform the integration over $N_{\bar\theta,0}$ first. Since
$\eta$ is trivial on $N_{\bar\theta,0}$ this integration is over
$\1_V$ only and it just gives a scalar factor. It thus remains to
discuss
$$
\int_{A_{CM}} a^{w_\theta\chi+\la_0+\la-\rho+2\rho_\theta}
\int_{\prod_{\alpha\in\Delta\setminus\bar\theta}N_\alpha}\eta(^an^{-1})\1_{n_0
V}(n) dnda.
$$
We may assume that
$n_0\in\prod_{\alpha\in\Delta\setminus\bar\theta}N_\alpha$ and
that $V$ is a product
$V=V_1\times\prod_{\alpha\in\Delta\setminus\bar\theta} V_\alpha$.
Splitting all integrals into respective products over
$\Delta\setminus\bar\theta$ we get that the above equals
$$
\prod_{\alpha\in\Delta\setminus\bar\theta}\left(
\int_{A_\alpha}a^{w_\theta\chi+\la_0+\la-\rho+2\rho_\theta}\int_{N_\alpha}\eta(^an^{-1})\1_{n_0V_\alpha}(n)
dn da\right).
$$
Fix $\alpha\in\Delta\setminus\bar\theta$ and consider the factor
attached to $\alpha$. Since $\eta$ is generic this factor equals
$$
\int_{a_\alpha'A_\alpha^-}a^{w_\theta\chi+\la_0+\la-\rho+2\rho_\theta}\int_{V_\alpha}\eta(^an^{-1})dn
\eta(^an_0^{-1}) da,
$$
for some $a_\alpha'\in A_\alpha$. There is a second
$a_\alpha''\in A_\alpha$ with $a_\alpha'' A_\alpha^-\subset
a_\alpha'A_\alpha^-$ and $\eta(^an^{-1})=1=\eta(^an_0^{-1})$ for
every $n\in V_\alpha$. Thus the contribution under consideration
is a sum of two summands the first of which is an integral
$$
\int_D
a^{w_\theta\chi+\la_0+\la-\rho+2\rho_\theta}\int_{V_\alpha}\eta(^an^{-1})dn\eta(^an_0^{-1})
da,
$$
where $D\subset A_\alpha$ is compact open and the function
$$
a\ \mapsto\ \int_{V_\alpha}\eta(^an^{-1})dn\eta(^an_0^{-1})
$$
is locally constant. This summand thus turns out to be a
polynomial in $q^{\pm\la_\alpha}$. The second summand is a
constant times
$$
\int_{a_\alpha''
A_\alpha^-}a^{w_\theta\chi+\la_0+\la-\rho+2\rho_\theta}da,
$$
which is rational in $q^{\la_\alpha}$. This finishes the proof of
the theorem.
\qed

\begin{conjecture}
We conjecture that the Mellin transform $I_{\la_0,f}(\la)$
always is a rational function in the $q^{\la_\alpha}$.
\end{conjecture}

We now recall the notion of an unramified group. The group $G$ is
called {\it unramified} if the two equivalent conditions hold:
\begin{itemize}
\item $G$ is obtained by base change from a smooth reductive group
scheme over ${\rm Spec}\CO$, where $\CO$ is the ring of integers
of $F$.
\item
$G$ is quasisplit over $F$ and split over an unramified extension.
\end{itemize}

\begin{proposition}
The conjecture holds if $G$ is unramified and of $F$-rank one,
and $\pi=\pi_\nu$ for some unramified quasicharacter $\nu$ of $A$.
\end{proposition}

\prf For $G$ and $\pi$ unramified Casselman and Shalika
\cite{CassShal} gave an explicit expression for $W_{f_0}(a)$,
where $f_0$ is the standard class one vector in $V_\pi$. This
formula indicates that the conjecture holds for $f=f_0$. For
arbitrary $f$ we have that $f-f(1)f_0$ vanishes at $1$ and
therefore lies in $V_\pi^c$, since $G$ is of rank one. By Theorem
\ref{main} and linearity the proposition follows.
\qed

\section{The group $\GL_2(F)$}
In this section we will compute the Mellin transform
$I_v(\la)$ for a vector $v$ which is fixed under the
maximal compact subgroup $K$. Since the Mellin transform
is taken over the derived torus it suffices to do the
computations in the group $\SL_2(F)$. We therefore let $P$
be the parabolic subgroup of all upper triangular matrices
in $\SL_2(F)$. We let $\pi=\pi_\nu$, $\nu\in\C$, the
unramified principal series representation. The space of
$\pi$ consists of all locally constant functions
$f:\SL_2(F)\ra\C$ satisfying $ f(a(t)nk)\= |t|^{\nu+1}f(k),
$ where $a(t)=\matrix t{}{}{\rez t}$ for $t\in F^\times$
and $k\in \SL_2(\CO)$. In particular let $f_0$ be the
standard class 1 vector, so $ f_0(a(t) nk)\= |t|^{\nu+1}. $
Then, if $\Re(\nu)<0$, we are in the range of Theorem
\ref{main}. We let
$$
\psi(f)\= \int_N \eta(n^{-1})f(wn) dn.
$$
We want to compute
$$
W_{f_0}(a(t))\=\int_N\eta(n^{-1})f_0(wna(t)) dn.
$$
We need to make the $AN K$-decomposition of $\SL_2(F)$ more
explicit. So assume that $g\in G$ is written as
$$
g\=\matrix abcd\=\matrix{t}{}{}{t^{-1}} \matrix 1x01 \matrix
\alpha\beta\gamma\delta,
$$
where $x\in F$ is arbitrary and
$\matrix\alpha\beta\gamma\delta\in \SL_2(\CO)$. Then
$$
g\=\matrix abcd\= \matrix{*}{*}{t^{-1}\gamma}{t^{-1}\delta}.
$$
So the fractional ideal generated by $c$ and $d$ is $t^{-1}\CO$.
For $n\=\matrix 1x{}1\in N$ this means
$$
f_0(wn)\=\begin{cases} 1 & {\rm if}\ |x|\le 1,\\ |x|^{-\nu-1} &
{\rm if}\ |x|>1.\end{cases}
$$
This leads to
$$
f_0(na(t))\=\begin{cases}|t|^{-\nu-1} & {\rm if}\ x\in t^2\CO,\\
|t|^{1+\nu}|x|^{-\nu-1} & {\rm if\ } x\ni t^2\CO.\end{cases}
$$
Now assume that $\eta\left(\matrix 1x{}1\right)=\tilde\eta(x)$
for some character $\tilde\eta$ of $F$ which is trivial on $\CO$
but nontrivial on $\varpi^{-1}\CO$. Then we get
\eqn
W_{f_0}(a(t)) &=& \int_N \eta(n^{-1})f_0(Wna(t)) dn\\
&=& |t|^{-\nu-1}\int_{t^2\CO}\tilde\eta(-x) dx\\
&& + |t|^{\nu+1}\sum_{j=1}^\infty
q^{j(-\nu-1)}\int_{\varpi^{-j}t^2\CO^\times}\tilde\eta(-x) dx
\endeqn
Write
\eqn
\int_{\varpi^{-j}t^2\CO^\times}\tilde\eta(-x) dx &=&
\int_{\varpi^{-j}t^2\CO}\tilde\eta(-x) dx -
\int_{\varpi^{1-j}t^2\CO}\tilde\eta(-x) dx\\
&=& |\varpi^{-j}t^2|\1_\CO(\varpi^{-j}t^2) -
|\varpi^{1-j}t^2|\1_\CO(\varpi^{1-j}t^2)
\endeqn
to get
\eqn
W(a(t))&=& |t|^{-\nu+1} \1_\CO(t) + |t|^{-\nu+1}\sum_{j=1}^\infty
q^{-\nu j}(\1_{\varpi^j\CO}(t^2)-\rez
q\1_{\varpi^{j-1}\CO}(t^2))\\
&=&
\frac{1-q^{-\nu-1}}{1-q^{-\nu}}(|t|^{1-\nu}-q^{-\nu}|t|^{1+\nu})
\endeqn
This implies that with
$$
I(s)\= \int_{F^\times}W(a(t))|t|^{s-1} d^\times t
$$
we get
$$
I(s)\=
\frac{(1-q^{-\nu-1})(1+q^{-s})}{(1-q^{-(s+\nu)})(1-q^{-(s-\nu)})}.
$$
We want to compare this with automorphic $L$-functions. To
choose an $L$-factor for $\bar\pi_\nu$ we have to fix a
representation of the Langlands dual group of $\SL_2(F)$
which is $\PGL_2(\C)$. The smallest dimensional
representations of this group are the trivial
representation and the adjoint representation $\Ad :
\PGL_2(\C)\ra \GL_3(\C)$. the $L$-factors are:
$$
L(\pi_\nu,triv,s)\= \rez{1-q^{-s}}
$$
and
$$
L(\pi_\nu,\Ad,s)\=
\rez{(1-q^{-(s+\nu)})(1-q^{-s})(1-q^{-(s-\nu)})}.
$$
We have proved the following proposition.

\begin{proposition} For $G=\GL_2(F)$ the Mellin transform of the
unramified Whittaker function equals
$$
I(s)\=\frac{L(\pi,\Ad,s)}{L(\pi,triv,s)}F_\nu(s),
$$
where $F_\nu(s)$ is the entire function
$F_\nu(s)=(1-q^{-\nu-1})(1+q^{-s})$.
\end{proposition}

We compare our Mellin transform to the one considered by Jacquet
and Langlands in \cite{JL}. They prove (Theorem 2.18), that the
integral defining
$$
\Psi(s,v)\=\int_{F^\times}W_v(\matrix t{}{}1 |t|^{s-\rez
2}d^\times t
$$
converges if the real part of $s$ is large enough. Further there
is a unique Euler factor $L(\pi,s)$, depending on $\pi$ and not
on $v$, such that
$$
\Psi(s,v)\= L(\pi,s)\Phi(s,v),
$$
where $\Phi(s,v)$ is entire in $s$ and there is a $v_0\in V_\pi$
such that $\Phi(s,v_0)=a^s$ for some positive constant $a$.

To compare this with $I_v(s)$ let $c=\#
(\CO^\times)/(\CO^\times)^2$. Then
$$
\Psi(s,v)\=\rez c \sum_{\ga:F^\times /(F^\times)^2}
\int_{F^\times} W_{\pi(\matrix \ga {}{}1 v}\left(\matrix
{t^2}{}{}1\right) |t|^{2s-1} d^\times t.
$$

This gives the following Proposition.

\begin{proposition}
If the central torus of all multiples of the identity matrix acts
trivially on $V_\pi$ then we have
$$
\Psi(s,v)\= I_{F(v)}(2s).
$$
where $F(v)=\rez c\sum_{\ga:F^\times /(F^\times)^2}\pi(\matrix
\ga {}{}1 v$.
\end{proposition}

\section{The unramified transform for $\GL_n(F)$}

For the case of $\GL_n(F)$ for $n>2$ one could either do
similar explicit computations or use the main result of
\cite{CassShal} to derive the following result:

\begin{proposition}
Assume $\pi=\pi_\nu=\Ind_P^G(1\otimes\nu\otimes 1)$, then
\eqn
I(\la)&=&
\left(\prod_{\alpha\in\Delta}\frac{1-q^{-\nu_\alpha-1}}{1-q^{\nu_\alpha}}\right)\\
&&\times\sum_{w\in W} \sign(w)
\left(\prod_{^{\alpha\in\Delta}_{w\alpha
<0}}q^{\nu_\alpha}\right)\left(
\prod_{\alpha\in\Delta}\rez{1-q^{-(\la_\alpha+(w\nu)_\alpha)}}\right)
\endeqn
\end{proposition}

We consider the special values at $\la_\alpha=s\in\C$ for all
$\alpha\in\Delta$. Write the corresponding function as
$I(s)=I(s,\dots,s)$. Further let
$$
C(\nu)\=\prod_{\alpha\in\Delta}\frac{1-q^{-\nu_\alpha-1}}{1-q^{\nu_\alpha}}
$$
and
$$
D(\nu,w)\=\prod_{^{\alpha\in\Delta}_{w\alpha
<0}}q^{\nu_\alpha},
$$
then
$$
I(s)\= C(\nu)\sum_{w\in
W}\sign(w)D(\nu,w)\prod_{\alpha\in\Delta}\rez{1-q^{-(s+(w\nu)_\alpha)}}.
$$
We get
\eqn
I(s)&=& C(\nu) \prod_{\beta\in\Phi}\rez{1-q^{-(s+\nu_{\beta})}}\\
&&\times\sum_{w\in W}\sign(w) D(\nu,w)\prod_{\beta\in\Phi\setminus
w\Delta}(1-q^{-(s+\nu_\beta)}).
\endeqn
Set
$$
F_\nu(s)\= C(\nu) \sum_{w\in W}\sign(w)
D(\nu,w)\prod_{\beta\in\Phi\setminus
w\Delta}(1-q^{-(s+\nu_\beta)})
$$
and
$$
L_\pi(s)\=\frac{L(\pi,\Ad,s)}{L(\pi,triv,s)^{n-1}},
$$
then we have

\begin{proposition}
The function $F_\nu(s)$ is entire and
$$
I(s)\= L_\pi(s) F_\nu(s).
$$
\end{proposition}

\section{The general unramified case}
Let $G$ be unramified. For $\alpha\in\Delta$ consider the
next to minimal parabolic $P_{\{\alpha\} }=M_{\{\alpha\}
}N_{\{\alpha\} }$. Then $M_{\{\alpha\} }$ has semisimple
rank one. So the simply connected covering $\tilde
M_{\{\alpha\} }$ of its derived group has rank one. There
are only two possible cases (see \cite{CassShal}): either
$\tilde M_{\{\alpha\} }\cong \SL_2(F)$ or $\tilde
M_{\{\alpha\} }\cong SU_3(F)$. For $\alpha\in\Delta$ let
$$
\eps_\alpha\=\begin{cases} 0 & {\rm if}\ \tilde
M_{\{\alpha\} }\cong \SL_2(F),\\ 1 & {\rm if}\ \tilde
M_{\{\alpha\} }\cong \SU_3(F).\end{cases}
$$
For $\nu\in\a^*$ let
$$
C(\nu)\=\prod_{\alpha\in\Delta}
\frac{(1-(-1)^{\eps_\alpha}q^{-\nu_\alpha-1})(1-q^{-\nu_\alpha-2})^{\eps_\alpha}}
     {1-q^{(1+\eps_\alpha)\nu_\alpha}}.
$$
For $w\in W$ let
$$
D(\nu,w)\=\prod_{\frack{\alpha\in\Delta}{w\alpha
<0}}q^{(1+\eps_\alpha)\nu_\alpha}.
$$

\begin{theorem}
For an unramified representation $\pi_\nu$, $\nu\in\a^*$
the Mellin transform of the Whittaker function attached to
the class one vector equals
$$
I(\la)\= C(\nu)\sum_{w\in W}\sign(w)
D(\nu,w)\prod_{\alpha\in\Delta}\rez{1-q^{-(\la_\alpha+(w\nu)_\alpha)}}.
$$
In particular, for $I(s)=I(s,\dots,s)$ we get
$$
I(s)\= L_\pi(s) F_\nu(s),
$$
where $F_\nu(s)$ is the entire function
$$
F_\nu(s)\= C(\nu)\sum_{w\in W}\sign(w)
D(\nu,w)\prod_{\beta\in \phi\setminus
w\Delta}(1-q^{-(s+\nu_\alpha)}),
$$
and $L_\pi(s)$ is the Euler factor
$$
L_\pi(s)\=\frac{L(\pi,\Ad,s)}{L(\pi, triv,s)^{\dim A'}}.
$$
\end{theorem}

{\small
University of Exeter, Mathematics, Exeter EX4 4QE, Devon, UK\\
a.h.j.deitmar@ex.ac.uk}

\end{document}